\documentclass[11pt]{article}

\usepackage{amssymb,a4wide}
\usepackage{url}

\newcommand{\ndv}{v}

\newcommand{\schm}{\mathsf{C}} 
\newcommand{\schmd}{\mathsf{D}}

\newcommand{\bridge}[1]{\overbracket[0.6pt]{\,#1\,}} 
\newcommand{\bridget}{\bridge{\ndt}}

\newcommand{\bridgetree}{\bridge{\tree}}

\newcommand{\sidebyside}[2]{
    \hbox to \textwidth{
                                            \vtop{\hsize=.45 \textwidth \parindent=0pt \centering #1 \vskip1sp}
                                            \hfill
                                            \vtop{\hsize=.45 \textwidth \parindent=0pt \centering #2}
                                            } }

\newcommand{\defstyle}[1]{\textup{\textbf{#1}}}

\newcommand{\cmp}[1]{\textnormal{\textsf{(C#1)}}}

\newcommand{\tree}{\ensuremath{\mathfrak{T}}}
\newcommand{\trt}{\mathfrak{T}}
\newcommand{\trs}{\mathfrak{S}}

\newcommand{\forestf}{\mathfrak{F}}

\newcommand{\forrom}[2]{F\left( #1 \backslash \mathsf{#2} \right)}
\newcommand{\forroml}[1]{F_l{\left(#1\right)}}
\newcommand{\forromu}[2]{F_u{\left( #1 \backslash \mathsf{#2} \right)}}
\newcommand{\forfra}[2]{\mathfrak{F}\left( #1 \backslash \mathsf{#2} \right)}
\newcommand{\forfral}[1]{\mathfrak{F}_l{\left(#1\right)}}
\newcommand{\forfrau}[2]{\mathfrak{F}_u{\left( #1 \backslash \mathsf{#2} \right)}}

\newcommand{\treet}{(T; <)}

\newcommand{\nds}{s}
\newcommand{\ndt}{t}
\newcommand{\ndu}{u}

\newcommand{\tlx}[2]{{#1}^{< #2}}
\newcommand{\tleq}[2]{{#1}^{\leqslant #2}}
\newcommand{\tg}[2]{{#1}^{> #2}}
\newcommand{\tgeq}[2]{{#1}^{\geqslant #2}}

\newcommand{\pathp}{\mathsf{P}}
\newcommand{\patha}{\mathsf{A}}
\newcommand{\pathb}{\mathsf{B}}
\newcommand{\pathc}{\mathsf{C}}

\newcommand{\pathx}{\mathsf{X}}
\newcommand{\pathq}{\mathsf{Q}}

\newcommand{\segm}{\mathsf{I}}

\newcommand{\segmj}{\mathsf{J}}

\newcommand{\stem}{\mathsf{S}}
\newcommand{\stemr}{\mathsf{R}}

\newcommand{\setpaths}{\mathcal{P}}

\newcommand{\lina}{\mathfrak{A}}

\newcommand{\ndb}{b}

\newcommand{\linacl}{\mathcal{A}}
\newcommand{\linbcl}{\mathcal{B}}

\renewcommand{\epsilon}{\varepsilon}

\renewcommand{\setpaths}{\mathsf{Paths}}

\newtheorem{theorem}{Theorem}

\newtheorem{corollary}[theorem]{Corollary}

\newtheorem{definition}[theorem]{Definition}

\newtheorem{lemma}[theorem]{Lemma}

\newtheorem{proposition}[theorem]{Proposition}

\newtheorem{fact}[theorem]{Fact}
\newenvironment{proof}[1][Proof]{\textbf{#1.} }{\ \rule{0.5em}{0.5em}}

\usepackage{mathtools}
\usepackage[font={small}]{caption}
\usepackage{xcolor}

\title
{Structural theory of trees \\
II. Completeness and completions of trees}

\author{Valentin Goranko$^1$, Ruaan Kellerman$^2$, and Alberto Zanardo$^3$  \\ 
$^1$Stockholm University, \ $^2$University of Pretoria, \ $^3$University of Padova \\ 
\textrm{valentin.goranko@philosophy.su.se}, \ 
\textrm{ruaan.kellerman@up.ac.za}, \  
\textrm{alberto.zanardo@unipd.it}} 

\begin{document}

\maketitle

\begin{abstract}
Trees are partial orderings where every element has a linearly ordered set of smaller elements.  
We define and study several natural notions of completeness of trees, extending Dedekind completeness of linear orders and  Dedekind-MacNeille completions of partial orders. We then define constructions of \emph{tree completions} that extend any tree to a minimal one satisfying the respective completeness property.
\end{abstract}

\section{Introduction}
\label{sec:intro}

Trees are connected partial orderings where every element has a linearly ordered set of predecessors (smaller elements). They are ubiquitous structures, naturally arising in a wide variety of contexts in mathematics, computer science, game and decision theory, linguistics, philosophy, etc.; some applications in these fields are mentioned in  \cite{GorankoKellermanZanardo2021a}.

This work is a contribution to the general structural theory of trees, and a follow-up to 
 \cite{GorankoKellermanZanardo2021a} which initiated that study.
We emphasize that in our definition and study of trees we do not assume well-foundedness, which is a standard assumption in the prevailing tradition of set-theoretic studies of trees, generalizing and extending the theory of ordinals, cf.~e.g.~\cite{Jech}, \cite{Todorcevic}. We only point out here that the well-foundedness assumption makes a very substantial difference, both in the general theory and in the particular properties of trees, and without that assumption, the study of trees remains mostly order-theoretic and extends in a quite non-trivial way the theory of linear orderings comprehensively explored in \cite{Rosenstein}.  
Furthermore, whereas our general theory covers both finite and infinite trees, the notions of completeness that are considered in this paper become nontrivial only in the case of infinite trees. 

In the present paper we define and explore several natural notions of \emph{completeness of trees}, 
intuitively stating that there are no `gaps' or `missing nodes', in one or another sense. 
These notions of completeness of trees naturally extend Dedekind completeness of linear orders and Dedekind-MacNeille completions of partial orders which result in complete lattices. 
These are important and well-studied constructions, but there is limited literature on their application specifically to trees, where there are some essential subtleties. Dedekind completeness of linear orders can be equivalently defined in terms of the existence of suprema of all non-empty sets that are bounded above, and in terms of the existence of infima of all non-empty sets that are bounded below. However, in the case of trees these two characterisations differ substantially. Indeed, in a tree, infima of linearly ordered sets of nodes that are bounded below, are unique whenever they exist, while linearly ordered set of nodes that are bounded above may have suprema on some paths in the tree and not on others, and thus it may have several (if any) suprema in that tree. This leads to a variety of notions of completeness of trees, generally of the following type: consider a (natural) family $\mathcal{F}$ of sets of nodes in the tree, and say that the tree is \emph{$\mathcal{F}$-complete} if every set of nodes in tree that is bounded below and belongs to $\mathcal{F}$ has an infimum. This generic notion, applied to the family of all subsets, defines a standard generalisation of Dedekind completeness of trees, variations of which have been defined and studied e.g. in \cite{Droste85}, \cite{Rubin93},  \cite{Warren}, and \cite{Barham}. Moreover, the generic notion of $\mathcal{F}$-completeness also makes very good sense when applied to other important families of nodes, such as paths (leading to `pathwise completeness'), antichains (leading to `antichain completeness'), pairs of nodes (leading to `branching completeness'), etc. 

\smallskip
\textbf{Contributions of the paper.}
Here we study the notions of absolute and relativised completeness that are mentioned above and that are, in our view, the most natural and important. We show that they are generally different, yet related to each other, and that they can be characterised in a fairly uniform way in terms of suprema and infima of downward closed and upward closed parts of paths.   
We then define, in a relatively uniform way, and study, several generic constructions of \emph{tree completions}, that extend any tree to a minimal one satisfying the respective completeness property. Each of these transforms in a `minimal canonical way' any tree into a complete tree in the respective sense. In particular, we present alternative constructions for producing trees that are equivalent to the  Dedekind-MacNeille completion of \cite{Warren} (applied there to the wider class of `cycle-free partial orders') and the `ramification completion' of \cite{Barham}. 

We note that completions of trees have various applications to the general theory of trees, which go beyond the scope of the present paper, but we only mention here that both are used when axiomatising the first-order theories, and other logical theories, of some important classes of trees, and when proving the completeness of such axiomatisations, cf. \cite{GorankoKellerman2021}.

\smallskip
\textbf{Structure of the paper.}
First, we provide the necessary terminology and notation in Section \ref{sec:prelim}. Then 
we define, compare and characterise several notions of completeness in Sections \ref{sec:completeness} and \ref{sec:characterisations}.  In Section \ref{sec:tree-completions} we introduce several constructions of tree completions that correspond to the respective completeness properties defined in Section \ref{sec:completeness}.  We end with concluding remarks and suggestions for further study in Section \ref{sec:concluding}.

\section{Preliminaries}
\label{sec:prelim}

We define here some basic notions on trees, in order to fix notation and terminology. The reader may also consult \cite{Jech}, \cite{KellermanThesis}, and \cite{Kellerman2018} for further details.

The order types of the linear orders $\left(\mathbb{N};<\right)$, $\left(\mathbb{Z};<\right)$, $\left(\mathbb{Q};<\right)$ and $\left(\mathbb{R};<\right)$, where in each instance $<$ denotes the usual ordering of that set, will be denoted as $\omega$, $\zeta$, $\eta$ and $\lambda$ respectively.

An ordered set $\left(A;<\right)$, {with a strict partial ordering $<$},  
 is \defstyle{downward-linear} if for every $x \in A$ the set $\{ y \in A : y < x \}$ is linear;  it 
is \defstyle{downward-connected} if, for every $x,y \in A$, there exists $z \in A$ such that $z \leqslant x$ and $z \leqslant y$.  A \defstyle{forest} is a downward-linear partial order.  A \defstyle{tree} is a downward-connected forest\footnote{Note that we do not assume well-foundedness of trees, nor even existence of a root.}. 
A  \defstyle{subtree}   
of a forest 
$\forestf = \left(F;<\right)$ 
is any substructure $\tree = \left(T;<^T\right)$ of $\forestf$ which is a tree, i.e.,  where $T$ is a non-empty downward-connected subset of $F$ and $<^T$ is the restriction of $<$ to  $T$.

The elements of a tree $\left(T;<\right)$  are called \defstyle{nodes} or \defstyle{points}. 
If a tree has a $<$-minimal node, then it is unique (by downward-connectedness)  and is called the \defstyle{root} of the tree. The $<$-maximal nodes (if any exist) are called \defstyle{leaves} of the tree. 

We will define various notions and notation in terms of an arbitrarily fixed tree $\tree = \treet$. 
First, for any nodes $\ndt, \ndu \in T$ we define $\ndt \smile \ndu$ to mean that $\ndt < \ndu$ or $\ndt = \ndu$ or $\ndu < \ndt$. If this holds, we say that $\ndt$ and $\ndu$ are \defstyle{comparable} nodes.
If $\ndt < \ndu$, 
the intervals $(\ndt,\ndu)$, $(\ndt,\ndu]$, $[\ndt,\ndu)$ and $[\ndt,\ndu]$ are defined as usual. 
For instance, if $\ndt < \ndu$ then $(\ndt,\ndu] := \{ x \in T : \ndt < x \leqslant u \}$, etc. 
 We also define the sets $\tlx{T}{\ndt} := \{ x \in T : x < \ndt\}$, $\tleq{T}{\ndt} := \{ x \in T : x \leqslant \ndt \}$, $\tg{T}{\ndt} := \{ x \in T : \ndt < x \}$ and $\tgeq{T}{\ndt} := \{ x \in T : \ndt \leqslant x \}$. 
 We will use analogous notation for the respective substructures (as partial orderings) of the tree $\tree$ over these sets, for instance,  $\tlx{\tree}{\ndt}$ denotes $\left(\tlx{T}{\ndt};<\upharpoonright_{\tlx{T}{t}}\right)$, etc.

For  non-empty subsets $A, B \subseteq T$ we define  $A < B$ (resp. $A \leqslant B$, $A > B$, $A \geqslant B$) iff $x<y$ (resp. $x \leqslant y$, $x > y$, $x \geqslant y$) for all $x \in A$ and $y \in B$.  Instead of $\{x\} < B$ we will also write $x < B$, and similarly for other relations and singleton sets.
Then, we define the set $\tlx{T}{A} := \{ x \in T : x < A \}$ and likewise the sets $\tleq{T}{A}$, $\tg{T}{A}$ and $\tgeq{T}{A}$.  The substructures of $\trt$ that have these sets as their underlying sets will be denoted as $\trt^{<A}$, $\trt^{\leqslant A}$, $\trt^{>A}$ and $\trt^{\geqslant A}$ respectively.

More generally, given any subset $A$ of $T$, $\trt^A$ will denote the structure $\left(A;<\upharpoonright_{A}\right)$.

Note that, for any $A \not=\emptyset$, $\tlx{\tree}{A}$ and  $\tleq{\tree}{A}$ are linear orders and that $\tg{T}{A}$ and $\tgeq{T}{A}$ are empty when $A$ is not linearly ordered.  If $A$ is linearly ordered then $\tg{\tree}{A}$ and $\tgeq{\tree}{A}$, if non-empty, are forests, while for every node $\ndt$, $\tgeq{\tree}{\ndt}$ is a tree that is rooted at $\ndt$.

For ease of readability, the sets $T^{\leqslant A}$ and $T^{\geqslant A}$ of lower bounds and upper bounds of $A$ will sometimes be denoted as $L(A)$ and $U(A)$ respectively.  If $A = \left\{x_1,x_2,\ldots,x_k\right\}$ then $L\left(A\right)$ and $U\left(A\right)$ will be written simply as $L\left(x_1,x_2,\ldots,x_k\right)$ and $U\left(x_1,x_2,\ldots,x_k\right)$.

A linearly ordered set of nodes in a tree is called a \defstyle{chain}. 
A maximal chain is called a \defstyle{path}. 
A set of nodes $\patha$ is  \defstyle{downward-closed} if $z \in \patha$ whenever $y \in \patha$ and $z < y$; respectively, $\patha$ is  \defstyle{upward-closed} if $z \in \patha$ whenever $y \in \patha$ and $y < z$.  
A non-empty downward-closed linearly ordered set of nodes that is bounded above is called a \defstyle{stem}. 
A non-empty subset $\pathb$ of a path $\patha$ is called a \defstyle{branch} when it is bounded below and upward-closed within $\patha$ (i.e.~if $x \in \pathb$ and $y \in \patha$ with $x < y$ then $y \in \pathb$).  Clearly if $\patha$ is a path with $\patha = \pathb \cup \pathc$, where $\pathb$ and  $\pathc$ are disjoint, then $\pathb$ is a stem if and only if $\pathc$ is a branch.  Note that every two distinct paths in a tree intersect in a stem. 
 The set of paths containing the node $\ndt$ (resp. the stem $\stem$) will  be denoted by $\setpaths_{\ndt}$ (resp. $\setpaths_{\stem})$.  

A set of nodes $\patha$ is called \defstyle{convex} if $z \in \patha$ whenever $x,y \in \patha$ and $x < z < y$.  
A convex linearly ordered set of nodes is called a \defstyle{segment}. 
A \defstyle{bridge} is a non-empty segment  $\patha$ such that, for every path $\pathp$, either $\patha \subseteq \pathp$ or $\patha \cap \pathp$ is empty.  Note that every singleton set of nodes $\{ t \}$ is a bridge. 

For each node $\ndt$ in $\tree$, the maximal bridge in $\tree$ containing $\ndt$ will be denoted by $\bridget$. A tree in which all maximal bridges are singletons is called a \defstyle{condensed tree}. As shown in \cite{GorankoKellermanZanardo2021a}, the relation between nodes of belonging to the same maximal bridge is an equivalence relation, and the corresponding partition of the tree $\tree$ into a family of maximal bridges defines a condensed tree, $\bridgetree$, called the \defstyle{condensation} of $\tree$. Moreover, a tree is condensed if and only if it is isomorphic to its condensation. Further details on properties of tree condensations and condensed trees can be found in \cite{GorankoKellermanZanardo2021a}.

A set $X$ of nodes is an \defstyle{antichain} if $x \not\smile y$ for all distinct $x$ and $y$ in $X$. Note that the intersection of an antichain $X$ and a linearly ordered set of nodes $Y$ is either a singleton or the empty set. The second alternative is excluded when $X$ is a maximal (by inclusion) antichain and $Y$ is a path. 

A  \defstyle{lower (respectively, upper) bound} for a set of nodes $X$ in $\trt$ is an element $b$ such that $b \leqslant x$ (respectively, $b \geqslant x$) for every $x \in X$. 
An  \defstyle{infimum (respectively, supremum)} of $X$ is a greatest lower bound (respectively, least upper bound) of $X$. If it exists, the infimum (respectively, supremum) of $X$ is unique and will be denoted $\inf(X)$ (respectively, $\sup(X)$).

A linear order  $\left(W;<\right)$ is \defstyle{Dedekind complete} when every non-empty subset in $W$ that is bounded below, has an infimum.  Equivalently (cf.~\cite{Rosenstein}) $\left(W;<\right)$ is Dedekind complete when every non-empty subset in $W$ that is bounded above, has a supremum.  Examples of Dedekind complete linear orders are $\left(\mathbb{N};<\right)$, $\left(\mathbb{Z};<\right)$, and $\left(\mathbb{R};<\right)$; a non-example is $\left(\mathbb{Q};<\right)$.

\section{Notions of completeness of trees}
\label{sec:completeness}

\subsection{Completeness properties and inter-dependence results}

We define here several natural notions of completeness in a tree. Note that the two equivalent characterisations of Dedekind complete linear orders do not transfer as equivalent properties over trees, because a linearly ordered set of nodes in a tree which is bounded above may have suprema on some paths and not on others, and thus it may have several (or none at all) suprema in the tree. 
So a variety of completeness properties emerge here.

For some of the definitions we need the following terminology:

\begin{itemize}
\item a node $\ndt$ in a tree will be called a \defstyle{branching point} when $\ndt = \inf\{\ndu,\ndv\}$ for some incomparable nodes $\ndu$ and $\ndv$;   

\item given two distinct paths $\pathp$ and $\pathq$ in the tree, the supremum in $\pathp$ (or, in $\pathq$) of the set $\pathp \cap \pathq$, if it exists, will be called a \defstyle{weakly branching point}\footnote{Note that a weakly branching point may also be a true branching point, in sense of the previous definition.}. 
\end{itemize}

The examples of trees in Figures~\ref{Fig:Dependence1} -- \ref{Fig:Dependence3}  below will be used in the work for proving independence results. The tree in Fig.~\ref{Fig:Dependence1} consists of the linear order $\omega + \zeta$ with a single leaf node attached to each of its elements. Every node in $\omega + \zeta$ is then a branching point. In Fig.~\ref{Fig:Dependence2}, the tree  consists of a copy of $\omega$ with two disjoint copies of $\omega$ appended to it. The first elements of each of these two copies are weakly branching points, but they are not branching points. The tree in Fig.~\ref{Fig:Dependence3} consists of a copy of the non-positive rational numbers $\eta^{\leqslant 0}$, with two disjoint copies of the positive rational numbers $\eta^{>0}$ appended to it.  The number $0$ in $\eta^{\leqslant 0}$ is both a branching point and weakly branching point of the tree.

\begin{figure}[tb]
	
	\hspace{-1cm}
 \begin{minipage}{0.32\textwidth}
     \centering
    \begin{picture}(40,100)
	\multiput(0,0)(0,12){6}{\circle*{3}}
	\put(0,0){\line(0,1){70}}
	\multiput(0,0)(0,12){6}{\line(1,1){10}}
	\multiput(10,10)(0,12){6}{\circle*{3}}
	\multiput(0,77)(0,7){3}{\circle{2}}
	\put(-15,40){$\omega$}
	\multiput(0,108)(0,12){5}{\circle*{3}}
	\put(0,98){\line(0,1){68}}
	\multiput(0,108)(0,12){5}{\line(1,1){10}}
	\multiput(10,118)(0,12){5}{\circle*{3}}
	\multiput(0,173)(0,7){3}{\circle{2}}
	\put(-14,135){$\zeta$}
	\end{picture} \hspace{0.1cm}
	
\caption{\label{Fig:Dependence1}}
   \end{minipage}
   \begin{minipage}{0.32\textwidth}
    \centering
     \begin{picture}(50,100)
	\put(0,0){\circle*{3}}
	\put(0,0){\line(0,1){70}}
	\multiput(0,77)(0,7){3}{\circle{2}}
	\put(-12,45){$\omega$}
	\put(5,98){\circle*{3}}
	\put(5,98){\line(2,3){40}}
	\multiput(49,164)(4,6){3}{\circle{2}}
	\put(38,135){$\omega$}
	\put(-5,98){\circle*{3}}
	\put(-5,98){\line(-2,3){40}}
	\multiput(-49,164)(-4,6){3}{\circle{2}}
	\put(-46,135){$\omega$}
	\end{picture}
	
 \caption{\label{Fig:Dependence2}}
   \end{minipage}
	\hspace{0.1cm}
   \begin{minipage}{0.32\textwidth}
     \centering
   
	\hspace{1cm}
  \begin{picture}(0,100)
	\multiput(0,21)(0,-7){3}{\circle{2}}
	\put(0,28){\line(0,1){70}}
	\put(0,98){\circle*{3}}
	\put(4,55){$\eta^{\leqslant 0}$}
	\multiput(5,105)(4,6){3}{\circle{2}}
	\put(17,123){\line(2,3){30}}
	\multiput(51,174)(4,6){3}{\circle{2}}
	\put(36,138){$\eta^{> 0}$}
	\multiput(-5,105)(-4,6){3}{\circle{2}}
	\put(-17,123){\line(-2,3){30}}
	\multiput(-51,174)(-4,6){3}{\circle{2}}
	\put(-56,138){$\eta^{> 0}$}
	\end{picture} \hspace{1cm}

 \caption{\label{Fig:Dependence3}}
   \end{minipage}
	\hfill
\end{figure}

We call a tree:

\begin{enumerate}
	\item[\cmp{1}] \defstyle{(Dedekind) complete} when every non-empty set of nodes that is boun\-ded below, has an infimum.
	
	\noindent {Hereafter, we will usually omit `Dedekind' and will simply write \defstyle{complete}.} 
		
	\item[\cmp{2}] \defstyle{pathwise (Dedekind) complete} when for every path $\pathp$ and each set $X \subseteq \pathp$ that is bounded below, $X$ has an infimum in $\pathp$. 	 
	
	\noindent {Again, hereafter, we will usually simply write \defstyle{pathwise complete}.}  
	
	\noindent Note that \cmp{2} is equivalent to: 
	\item[\cmp{$2^{\prime}$}]
		for every path $\pathp$ and each set $X \subseteq \pathp$ that is bounded above, $X$ has a supremum in $\pathp$;		
	\item[\cmp{3}] \defstyle{antichain complete} when every non-empty antichain that is boun\-ded below has an infimum;
	
	\item[\cmp{4}] \defstyle{branching complete}\footnote{In \cite{ZanardoBarcellanReynolds}, branching complete trees are called \emph{jointed} trees, in \cite{Barham} they are called \emph{ramification complete} trees, and in \cite{CourcelleDelhomme} they are called \emph{join-trees}.} when every pair of incomparable nodes (equi\-valently, every finite antichain) has an infimum;
	
	\item[\cmp{5}] \defstyle{weakly branching complete} when for any two distinct paths $\pathp$ and $\pathq$, the set $\pathp \cap \pathq$ has suprema in both $\pathp$ and $\pathq$;
	
	\item[\cmp{6}] \defstyle{weakly branching point complete} when for each path $\pathp$ and each non-empty set $X$ 
	of weakly branching points in $\pathp$, if $X$ is bounded below (respectively, above) then $X$ has an infimum (respectively, supremum)\footnote{Note that the two conditions, for the existence of infima and of suprema, are not equivalent. Indeed, consider a tree similar to the one in Figure~\ref{Fig:Dependence1} but with the leaves removed from the terminal $\zeta$-part of the tree.  The weakly branching points of this tree are precisely the non-leaf nodes in the initial $\omega$-part of the tree.  Clearly each non-empty set of weakly branching points has an infimum, while the set of all weakly branching points is bounded above but has no supremum. } in $\pathp$. 
\end{enumerate}

Note that each of the properties \cmp{1}, \cmp{2}, \cmp{3} and \cmp{4} fits the following generic notion of completeness: 
given a family $\mathcal{F}$ of sets of nodes in a tree $\trt$, we say that $\trt$ is \defstyle{$\mathcal{F}$-complete} if every set of nodes in $\trt$ that is bounded below and belongs to $\mathcal{F}$, has an infimum.  In the case of Dedekind completeness, take $\mathcal{F}$ to be the collection of all sets of nodes in $\trt$; in the case of pathwise completeness, take $\mathcal{F}$ to be the set of all linearly ordered sets of nodes in $\trt$; in the case of antichain completeness, take $\mathcal{F}$ to be the set of all antichains in $\trt$; and in the case of branching completeness, take $\mathcal{F}$ to be the set of all pairs of incomparable nodes in $\trt$. 

\smallskip

Some simple observations: 

\begin{itemize}
\item Every branching point is also a weakly branching point, while in branch\-ing complete trees, it also holds that every weakly branching point is a branching point.  
\item 
In any tree, every branching point is the greatest element of a maximal bridge, and every weakly branching point is either the greatest or the least element of a maximal bridge. The converses of these statements do not hold. For example, consider the tree that resembles the one in Figure~\ref{Fig:Dependence1} but with a single node $\ndu$ inserted between the initial portion $\omega$ and the terminal portion $\zeta$.  This node $u$ forms a maximal bridge (i.e.~$\bridge{u} = \{u\}$) and $u$ is both the least and the greatest element of that maximal bridge but $u$ is neither a branching point nor a weakly branching point.

\item In any tree, if $\pathp$ is a path and $X$ is a non-empty subset of $\pathp$ with infimum $\ndu$ in $\pathp$, then $\ndu$ will be the infimum of $X$ in \textit{every} path $\pathq$ that contains $X$; in other words, infima of non-empty chains are not path specific. 

\item  Moreover, in any \textit{branching complete} tree, if $\pathp$ is a path and $X \subseteq \pathp$ is a non-empty set with supremum $\ndv$ in $\pathp$, then $\ndv$ will be the supremum of $X$ in \textit{every} path $\pathq$ that contains $X$; in other words, suprema of non-empty chains in a branching complete tree are not path specific. 
\end{itemize}

It follows from these observations that in \textit{branching complete} trees, \cmp{6} is equivalent to the following property:
\begin{enumerate}
	\item[\cmp{7}]
		every non-empty chain of branching points that is bounded below (respectively, bounded above) has an infimum (respectively, supremum).
\end{enumerate}

The following theorem summarises all dependence results that hold between the properties \cmp{1} to \cmp{7}.

\begin{theorem} \label{Thm:Dependences} 
The following implications hold between the properties \cmp{1} to \cmp{7} and all implications that hold between \cmp{1} to \cmp{7} follow from these by transitivity. 
\begin{enumerate}
\item  \cmp{1} implies each of  
\text{\cmp{2} to \cmp{7}}; 

\item \cmp{2} is equivalent to \cmp{$2^{\prime}$};  

\item   \cmp{2}  implies \cmp{5} to \cmp{7};

\item   \cmp{3}  implies  \cmp{4}; 

\item   \cmp{4} implies \cmp{5}.

\end{enumerate}
\end{theorem}

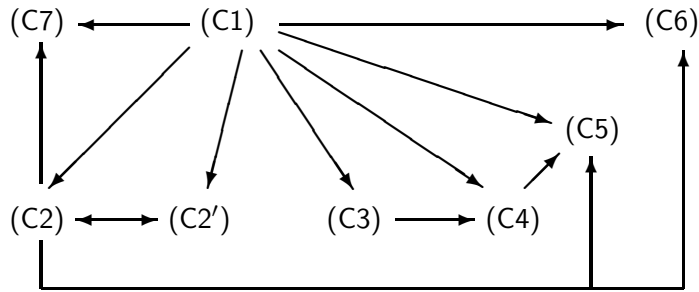
\begin{figure}[th]
	\begin{center}
		\begin{picture}(300,130)
			\put(72,100){\cmp{1}}
			\put(0,25){\cmp{2}}
			\put(60,25){\cmp{2$'$}}
			\put(120,25){\cmp{3}}
			\put(180,25){\cmp{4}}
			\put(210,59){\cmp{5}}
			\put(240,100){\cmp{6}}
			\put(0,100){\cmp{7}}
			\thicklines
			\put(68,102){\vector(-1,0){42}}
			\put(12,42){\vector(0,1){53}}
			\put(68,93){\vector(-1,-1){53}}
			\put(88,92){\vector(-1,-4){13}}
			\put(95,92){\vector(2,-3){35}}
			\put(102,92){\vector(3,-2){77}}
			\put(102,99){\vector(3,-1){103}}
			\put(102,102){\vector(1,0){130}}
			\put(195,40){\vector(1,1){13}}
			\put(26,28){\vector(1,0){30}}
			\put(56,28){\vector(-1,0){31}}
			\put(146,28){\vector(1,0){30}}
			\put(12,2){\line(0,1){18}}
			\put(12,2){\line(1,0){243}}
			\put(255,2){\vector(0,1){90}}
			\put(220,2){\vector(0,1){50}}
		\end{picture}
	\end{center}

	\caption{The dependencies that hold between \cmp{1} to \cmp{7}.\label{Fig:Dependencies}}
\end{figure}

\begin{proof}

We first consider the relationships between the properties \cmp{1} to \cmp{6}, after which we will relate them with the auxiliary property \cmp{7}.

1. Most implications from  \cmp{1} to the properties {\cmp{2} to \cmp{6}} are trivial.  For \cmp{5}, consider the sets of all upper bounds of $\pathp \cap \pathq$ in each of $\pathp$ and $\pathq$.

2. Just like the argument for linear orders: to show that \cmp{2} implies \cmp{$2^{\prime}$}, consider the set of all upper bounds of $X$ in $\pathp$ and take its infimum, and to show that \cmp{$2^{\prime}$} implies \cmp{2}, consider the set of all lower bounds of $X$ in $\pathp$ and take its supremum. 

3. The implication from \cmp{2} to \cmp{5} is like that of \cmp{1} to \cmp{5}; the implication to $\cmp{6}$ is trivial.

4. \cmp{3} implies \cmp{4} is trivial.

5. Take two points, one on $\pathp \setminus \pathq$ and the other on $\pathq \setminus \pathp$. Their infimum is a supremum of $\pathq \cap \pathp$ in each path.

\smallskip

Now, we show that all implications that do not follow from those listed in the proposition, do not hold.

First, consider the tree in Fig.~\ref{Fig:Dependence2}.  The properties \cmp{2}, and hence \cmp{5} and \cmp{6}, hold in this tree, while none of \cmp{1}, \cmp{3} and \cmp{4} hold in it.  It follows that none of the properties \cmp{2}, \cmp{5} and \cmp{6} imply any of \cmp{1}, \cmp{3} and \cmp{4}.

Next, consider the tree in Fig.~\ref{Fig:Dependence3}.  This tree shows that \cmp{3} (hence also \cmp{4}) does not imply \cmp{1}, and none of \cmp{3} to \cmp{6} imply \cmp{2}.

Similar to the previous example, consider the tree that consists of a copy of the rational numbers $\eta$, with two disjoint copies of the rational numbers $\eta$ appended to it.  This tree has the same structure as the one in Fig.~\ref{Fig:Dependence3} except that the weakly branching points of that tree are missing from it.  The property \cmp{6} holds vacuously in this tree, while \cmp{5} does not hold in it, hence \cmp{6} does not imply \cmp{5}.

To see that \cmp{4} (hence also \cmp{5}) does not imply \cmp{6}, consider the tree in Fig.~\ref{Fig:Dependence1}, which also shows that \cmp{4} does not imply  \cmp{3}.

The tree obtained from the one on Fig.~\ref{Fig:Dependence1} by removing the leaf nodes that were attached to the elements of the copy of $\zeta$ in $\omega + \zeta$, shows that \cmp{3} does not imply \cmp{6}.

This concludes the determination of all possible dependencies between the properties \cmp{1} to \cmp{6}.

\smallskip

Now, regarding the dependencies between \cmp{1} to \cmp{6} and  \cmp{7} on the class of \textit{all} trees, first note that \cmp{1} and each of the equivalent statements \cmp{2} and \cmp{2'}, implies \cmp{7}.  Further, none of \cmp{3}, \cmp{4}, \cmp{5} and \cmp{6} imply \cmp{7}.  To see that \cmp{3} does not imply \cmp{7}, consider the tree that consists of the path $\omega+\zeta$ with an extra node adjoined to the side of each node in the copy of $\omega$.  This tree is antichain complete but does not satisfy \cmp{7}.  That neither of \cmp{4} and \cmp{5} implies \cmp{7} can be seen from the tree in Figure~\ref{Fig:Dependence1}.   To see that \cmp{6} does not imply \cmp{7}, consider the tree obtained from the one in Figure~\ref{Fig:Dependence2} by adding a leaf to the side of each node in that tree.  Each node from the original tree is a branching node (hence also weakly branching node) in this new tree.  This tree satisfies \cmp{6} but not \cmp{7} as the set of branching nodes that lie on the initial copy of $\omega$, has no supremum.

In the other direction, \cmp{7} does not imply any of \cmp{1} to \cmp{6}.  Indeed, the tree in Figure~\ref{Fig:Dependence2} satisfies \cmp{7} but not \cmp{1}, \cmp{3} or \cmp{4}, while the tree in Figure~\ref{Fig:Dependence3} satisfies \cmp{7} but not \cmp{2}.  The tree that resembles the one in Figure~\ref{Fig:Dependence3} but with the greatest node in the initial copy of $\eta^{\leqslant 0}$ removed, vacuously satisfies \cmp{7}, but it does not satisfy \cmp{5}.  Finally, the tree that consists of a copy of $\omega \cdot (\zeta + \zeta)$ (i.e.~$\zeta + \zeta$ copies of $\omega$ placed end to end) with another copy of $\omega$ placed aside each copy of $\omega$ in $\omega \cdot (\zeta + \zeta)$, vacuously satisfies \cmp{7}, but does not satisfy \cmp{6}.
\end{proof}

\subsection{Rubin completeness}
\label{sec:Rubin}

A somewhat more restrictive notion of completeness is considered in \cite{Rubin93} (Def.~0.2). We say that a tree $\tree$ is \defstyle{Rubin complete} if
\begin{itemize}
\item[(1)] every non-empty set of nodes has an infimum; 
\item[(2)] every non-empty chain has a supremum;
\item[(3)]  if $\patha$ and $\pathb$ are disjoint non-empty convex chains, then $\inf(\patha) \not= \inf(\pathb)$. 
\end{itemize}
Then,  Rubin complete trees are (Dedekind) complete and hence they can be placed above \cmp{1} in the hierarchy considered in this section. It can also be observed that, by (1), Rubin complete trees are rooted and, by (2), every path in them has a leaf.

\begin{proposition} A complete tree is Rubin complete if and only if:   
\begin{itemize}
\item[(a)] 
for every set $\{ \pathp_i : i \in I\}$ of paths, if $\bigcap_{i \in I} \pathp_i \not= \emptyset$, then this set  has a maximum, 

\item[(b)]
 for any two paths $\pathp \not= \pathq$, $\pathp \setminus \pathq$ and $\pathq \setminus \pathp$ have a minimum, and 
 
 \item[(c)]
 every non-leaf node has an immediate successor on each path passing through it.
 \end{itemize}
\end{proposition}

\begin{proof} $(\Rightarrow)$ 
Let $\tree$ be a Rubin complete tree. Consider a set $\{ \pathp_i : i \in I\}$ of paths in $\tree$ and assume that  $\patha =  \bigcap_{i \in I} \pathp_i  \not= \emptyset$.  By condition (2), $\patha$ has a supremum $\ndt_\patha$. If $I$ is the singleton $\{i\}$, $\patha$ is $\pathp_i$ and $\ndt_\patha$ is the the leaf of it. If $|I| \geqslant 2$, assume for reductio  $\ndt_\patha \not\in \patha$, so that $\ndt_\patha \not\in \pathp_i$ for some $i \in I$. This implies $\ndt_\patha \not\leqslant \ndv$ for every $\ndv \in \pathp_i$. But $\pathp_i$ contains upper bounds of $\patha$.

Let now $\pathp \not= \pathq$ be paths in $\tree$. Consider the non-empty segments $\patha = \pathp \setminus \pathq$ and $\pathb = \pathq \setminus \pathp$ and let $\ndt_\patha$ and $\ndt_\pathb$ be their infima. Let $\ndt$ be the maximum of $\pathp \cap \pathq$, which exists by condition (a).
Assume for reductio that $\patha$ has no minimum, so that $\ndt_\patha = \ndt$. The node $\ndt$ is also the infimum of the segment $\pathb' = \{\ndt\} \cup  \pathq \setminus \pathp$ and hence $\inf(\patha) = \inf(\pathb') = \ndt$. But this contradicts (3) because  $\patha$ and $\pathb'$ are disjoint.

Finally, consider a non-leaf node $\ndt$ and a path $\pathp$ passing through it. Let $\patha$ be the set of all nodes $\ndu \in \pathp$ such that $\ndt < \ndu$. Then $\patha$ and $\pathb = \{ \ndt \}$ are disjoint non-empty convex chains. Since $\inf(\patha) \not= \inf(\pathb) = \ndt$, $\inf(\patha)$ must be the immediate successor of $\ndt$ in $\pathp$.

\smallskip
$(\Leftarrow)$ Let $\tree$ be a complete tree in which (a), (b), and (c) hold. Consider a non-empty chain $\patha$ in $\tree$ and the set $\{ \pathp_i : i \in I\}$ of all paths containing $\patha$. Call $\pathb$ the set  $\bigcap_{i \in I} \pathp_i$ and $\ndt_\pathb$ the maximum of $\pathb$. 
Given any $\pathp_i$, the pathwise completeness of $\tree$ implies that $\patha$ has a supremum $\ndu_i$ in $\pathp_i$. We can observe now that
 $\ndt_\pathb$ is an upper bound of $\patha$ and belongs to every $\pathp_i$. Then $\ndu_i \leqslant \ndt_\pathb$ for each $i$ and hence, for all $i, j$, $\ndu_i = \ndu_j$. This means that $\ndu_i$ is the supremum of $\patha$ also in $\tree$.

Consider now two disjoint segments $\patha$ and $\pathb$ and let $\ndt_\patha$ and $\ndt_\pathb$ be their infima. Two cases can be distinguished.
\\
Case 1: $\patha \not< \pathb$ and $\pathb \not< \patha$.  Then there exist two different paths $\pathp$ and $\pathq$ containing $\patha$ and $\pathb$, respectively. Since $\patha$ and $\pathb$ are disjoint, one of the inclusions $\patha \subseteq \pathp \setminus \pathq$ and $\pathb \subseteq \pathq \setminus \pathp$ must hold. Assume, w.l.o.g.,  the first one. By (b), we can consider the minimum $\ndt$ of $\pathp \setminus \pathq$. Then $\ndt \leqslant \ndt_\patha \in \pathp \setminus \pathq$. This implies $\ndt_\patha\not= \ndt_\pathb$ because $\ndt_\pathb \in \pathq$.
\\
Case 2:  $\patha < \pathb$ (the case $\pathb < \patha$ is similar). If $\patha$ is not a singleton set, then  $\ndt_\patha \not= \ndt_\pathb$ trivially holds. Else, assume $\patha = \{\ndu\}$  and let $\pathp$ be any path containing $\pathb$, so that $\ndu \in \pathp$. By (c), $\ndu$ has an immediate successor  $\ndu'$ in $\pathp$. Then the claim follows by the inequalities $\ndt_\patha \leqslant \ndu < \ndu' \leqslant \ndt_\pathb$.
\end{proof}

Observe that the strong condition (c) is needed only for dealing with the very particular case $\patha = \{\ndu\}< \pathb$.  This can be avoided, for instance, by assuming that $\patha$ and $\pathb$ contain incomparable nodes, or by assuming that chains are not singleton sets.   In these cases, condition (c) and the last part of the proof above can be dropped.

\subsection{Weakly branching completeness and partitions}
\label{sec:partitioning}

The property of weakly branching completeness, apart from being a natural notion of completeness in trees, also becomes critical when partitioning a tree along one of its stems.  Such partitions are used in \cite{GorankoKellerman2021} for the purpose of approximating trees as coloured linear orders\footnote{There the property of weakly branching completeness is simply called \textit{branching completeness} since no stronger form of branching completeness is used in that context.}.

 For $\tree$ any tree, $\patha$ any path in $\tree$ and $\nds \in \patha$, we define the following sets:
\begin{eqnarray*}
	\forroml{s} & := & \left\{ \begin{array}{cl}
	\emptyset, & \text{when $s$ has an immediate predecessor,} \\ ~ \\ 
	\displaystyle \left(\bigcap_{t<s}T^{>t}\right) \setminus T^{\geqslant s}, & \text{otherwise} \end{array} \right. \\
	\text{and} \\
	\forromu{s}{A} & := & T^{> s} \setminus \left(\,\bigcup_{t \in \patha \cap T^{>s}} \!\!\!\! T^{\geqslant t} \right).
\end{eqnarray*}
Provided that these sets are non-empty, we define the structures $\forfral{s} := \mathfrak{T}^{\forroml{s}}$ and $\forfrau{s}{A} := \mathfrak{T}^{\forromu{s}{A}}$.   In general, either or both of the sets $\forroml{s}$ and $\forromu{s}{A}$ may be empty, in which case the corresponding structure is left undefined.  The structures $\forfral{s}$ and $\forfrau{s}{A}$, if defined, are forests. They will be called respectively the \defstyle{lower side-forest} of $s$, and the \defstyle{upper side-forest} of $s$ with respect to $\patha$.  Finally, the \defstyle{side-forest} of $s$ with respect to $\patha$ is the forest $\forfra{s}{A} := \left(\mathfrak{T}^{\forrom{s}{A}};s\right)$ where
$$\forrom{s}{A} := \left\{s\right\} \cup \forroml{s} \cup \forromu{s}{A}.$$
The forests $\forfral{s}$ and $\forfrau{s}{A}$ are depicted in Fig.~\ref{Fig:Side-Forests}.

\begin{figure} 
	\begin{center}
	\begin{picture}(0,190)
	\put(0,-10){\line(0,1){195}}
	\put(0,80){\circle*{3}}
	\qbezier[60](5,85)(35,115)(65,145)
	\qbezier[60](5,85)(20,125)(35,165)
	\qbezier[25](65,145)(50,155)(35,165)
	\put(53,159){$\mathfrak{F}_{u}\left(s\backslash\mathsf{A}\right)$}
	\qbezier[60](5,75)(35,105)(65,135)
	\qbezier[60](5,75)(45,90)(85,105)
	\qbezier[25](65,135)(75,120)(85,105)
	\put(78,123){$\mathfrak{F}_{l}\left(s\right)$}
	\multiput(0,120)(0,10){3}{\circle*{3}}
	\multiput(-4,100)(0,5){3}{\circle{2}}
	\multiput(-4,150)(0,5){3}{\circle{2}}
	\multiput(0,40)(0,-10){3}{\circle*{3}}
	\multiput(-4,60)(0,-5){3}{\circle{2}}
	\multiput(-4,10)(0,-5){3}{\circle{2}}
	\put(-10,77){$s$}
	\put(-10,180){$\mathsf{A}$}
\end{picture}

	\caption{The side-forests $\forfral{s}$ and $\forfrau{s}{A}$.\label{Fig:Side-Forests}}
	\end{center}
\end{figure}

Intuitively, the upper side-forest of $s$ with respect to $\patha$ consists of those nodes that sit above $s$ but do not sit above any nodes on $\patha \cap T^{>s}$, and the lower side-forest of $s$ consists of those nodes that sit above $T^{<s}$ but are incomparable with $s$, unless $s$ has an immediate predecessor, in which case its lower side-forest is left empty so as not to coincide with the upper side-forest of that predecessor.

The following result is adapted from \cite{GorankoKellerman2021}.

\begin{proposition}
	Let $\trt$ be a weakly branching complete tree and let $\patha$ be a path in $\tree$.  The set $\displaystyle \left\{ \forrom{s}{A} \right\}_{s \in \patha}$ forms a partition of $T$.
\end{proposition}

\begin{proof}
	Each side-forest $\forrom{s}{A}$ is non-empty since $s \in \forrom{s}{A}$, and different side-forests $\forrom{s}{A}$ and $\forrom{t}{A}$ are disjoint by the way that side-forests are defined.  To see that $\displaystyle \left\{ \forrom{s}{A} \right\}_{s \in \patha}$ covers $T$, pick any node $u \in T$, and consider the following cases.
	
	Case 1: $u \in \patha$.  Then $u \in \forrom{u}{A}$.
	
	Case 2: $u \not\in \patha$.  Then there exists $v \in \patha$ such that $\left\{u,v\right\}$ is an antichain.  Let $\pathb$ be a path that contains $u$.  By the weakly branching completeness of $\tree$, the set $\patha \cap \pathb$ has an supremum $w$ in $\patha$.
	
	Case 2.1: $w < u$.  Then $u \in \forromu{w}{A}$.
	
	Case 2.2: $w \not< u$.  Then $u \in \forroml{w}$.
	
	\noindent This shows that $\bigcup_{s \in \patha} \forrom{s}{A} = T$.
\end{proof}

\section{Characterisations of the notions of completeness}
\label{sec:characterisations}

\begin{theorem} \label{Thm:DedekindCompleteness}
	A tree is complete \cmp{1} if and only if it is both pathwise complete \cmp{2} and branching complete \cmp{4}.
\end{theorem}

\begin{proof}
	($\Rightarrow$) \ Trivial.
	
	($\Leftarrow$) \ Suppose that \cmp{2}  and \cmp{4} hold in the tree $\trt$.  Consider a non-empty set $X$ of nodes in $\trt$ that is bounded below.  Let $\patha$ be any path in $\trt$ that intersects $X$.  Then $\patha \cap X$ is a non-empty subset of $\patha$ that is bounded below, hence, by \cmp{2},  it has an infimum on $\patha$, call it $\ndv$.
	For any path $\pathp$ that intersects with $X$, again, $\pathp \cap X$ has an infimum on $\pathp$, call it $\ndt_{\pathp}$.  By \cmp{4}, the set $\left\{\ndv,\ndt_{\pathp}\right\}$ has an infimum $\ndu_{\pathp}$, which is a node on $\patha$. Observe that $\ndt_\patha = \ndv = \ndu_\patha$.  Now consider the set
	\[ Y := \{ \ndu_{\pathp} \mid \, \pathp \mbox{ is any path that intersects with } X \}. \]
For any lower bound $\ndu$ of $X$, we have $\ndu \leqslant \ndv$, $\ndu \leqslant \ndt_{\pathp}$, and hence $\ndu \leqslant \ndu _{\pathp}$  for every path $\pathp$ intersecting $X$. 	
	Then, the set $Y$ is a non-empty subset of the path $\patha$ and is bounded below. By property \cmp{$2$} it has an infimum $\ndt$ in $\patha$.  We claim that $\ndt$ is the infimum of $X$. Indeed, the definition of $Y$ implies that $\ndt$ is a lower bound of $X$, and we have observed above that every lower bound $\ndt'$ of $X$ is also a lower bound of $Y$, and hence $\ndt' \leqslant \ndt$.
\end{proof}

\begin{corollary}
	A pathwise complete \cmp{2} tree $\trt$ is antichain complete \cmp{3} if and only if $\trt$ is branching complete \cmp{4}, if and only if $\trt$ is complete \cmp{1}.
\end{corollary}

The next result shows that the property of antichain completeness is equivalent to the conjunction of the properties of branching completeness and ``half'' of weakly branching point completeness (bearing in mind that in branching complete trees, branching points and weakly branching points coincide).

\begin{proposition} \label{Thm:Antichains}
	A tree $\trt$ is antichain complete \cmp{3} if and only if $\trt$ is branching complete \cmp{4} and every non-empty chain of branching points that is bounded below, has an infimum.
\end{proposition}

\begin{proof}
	($\Rightarrow$) \ Assume that every non-empty antichain in $\trt$ that is bounded below, has an infimum.  It is immediate that $\trt$ is branching complete.  Let $X$ be a non-empty chain of branching points in $\trt$ that is bounded below.  For each $x \in X$, let $\ndu_x$ and $\ndv_x$ be incomparable nodes such that $x = \inf\{\ndu_x,\ndv_x\}$ and let $Y := \bigcup_{x \in X} \{\ndu_x,\ndv_x\}$.  It follows from Zorn's Lemma 
	that there exists $Z \subseteq Y$ that is maximal with respect to being an antichain, and $Z$ must be bounded below since $X$ is bounded below.  Let $\ndt := \inf(Z)$.
	
	Observe that $L(X) = L(Y)$.  Let $\nds \in L(X)$.  Since $Z \subseteq Y$ then $L(Y) \subseteq L(Z)$ from which $\nds \in L(Z)$ hence $\nds \leqslant \ndt$.
	To see that $\ndt \in L(X)$, first note that, by the maximality of $Z$, there exists, for each $y \in Y$, a node $y' \in Z$ such that $y' \smile y$.  Let $x \in X$.  If $(\ndu_x)' \neq (\ndv_x)'$ then $(\ndu_x)' \not\smile (\ndv_x)'$ (since $Z$ is an antichain) hence $\ndt \leqslant \inf\{(\ndu_x)',(\ndv_x)'\} = \inf\{\ndu_x,\ndv_x\} = x$.  If, on the other hand, $(\ndu_x)' = (\ndv_x)'$, then $(\ndu_x)' \smile \ndu_x,\ndv_x$, hence $\ndt \leqslant (\ndu_x)' \leqslant \inf\{\ndu_x,\ndv_x\} = x$, hence $\ndt \in L(X)$, from which it follows that $t = \inf(X)$, as required.
	
	($\Leftarrow$) \ Let $X$ be a non-empty antichain that is bounded below.  Let $x_0 \in X$ and let $\pathp$ be a path that passes through $x_0$.  For each $y \in X \setminus \{x_0\}$, let $\ndu_y := \inf\{x_0,y\} \in \pathp$.  The set $\left\{ \ndu_y : y \in X \setminus \{x_0\} \right\}$ is then a chain of branching points that is bounded below, hence this set has an infimum $\ndv$.  It is readily verified that $\ndv$ is also the infimum of $X$.
\end{proof}

Completeness issues lead us to consider particular stems and branches.  
Recall that $\bridge{\ndt}$ denotes the maximal bridge in the tree containing $\ndt$.
We call a stem $\patha$ a \defstyle{trunk} when $\patha = \bigcup_{t \in \patha} \bridge{\ndt}$ 
and, similarly, call a branch $\pathb$ a \defstyle{limb} when $\pathb = \bigcup_{t \in \pathb} \bridge{\ndt}$.  For a non-example of a trunk and a limb, consider the tree $\trt$ in Figure~\ref{Fig:Dependence2}.  For each node $\ndu$ in the initial copy of $\omega$ in $\trt$, the set $T^{\leqslant \ndu}$ forms a stem that is not a trunk, and for each node $\ndv$ in each of the two copies of $\omega$ that are appended atop the initial copy of $\omega$, the set $T^{> \ndv}$ forms a branch that is not a limb.

If $\pathp$ is a path and $\patha \subseteq \pathp$ then $\patha$ is a trunk if and only if $\pathp \setminus \patha$ is a limb.  If $\patha$ is a trunk and $\pathb$ is a limb for which $\patha \cap \pathb = \emptyset$ and $\patha \cup \pathb$ is a path, then $\patha$ and $\pathb$ will be called \defstyle{complementary}.  In condensed trees, stems and trunks coincide, as do branches and limbs.  Note that if a trunk $\patha$ and limb $\pathb$ are complementary, then the supremum of $\patha$, if it exists, will be the infimum of $\pathb$.  However, if $\pathb$ has an infimum, that infimum need not be the supremum of $\patha$, in fact, $\patha$ need not even have a supremum. The tree in Figure~\ref{Fig:Dependence2} can serve as a counterexample.

\begin{theorem} \label{Thm:TrunksForm}
	Let $\trt$ be any tree.
	\begin{enumerate}
		\item
			For any non-singleton set $\left\{\pathp_i\right\}_{i \in I}$ of paths in $\trt$, the set $\displaystyle \bigcap_{i \in I} \pathp_i$, if non-empty, is a trunk.
		\item
			For any path $\pathq$ and any non-singleton set $\left\{\pathp_i\right\}_{i \in I}$ of paths in $\trt$, the set $\displaystyle \pathb := \pathq \setminus \bigcup_{i \in I} \pathp_i$, if non-empty, is a limb (hence $\pathq \setminus \pathb$ is a trunk).
		\item
		Every trunk in $\trt$ is of one, and possibly both, of the two forms above.
	\end{enumerate}
\end{theorem}

\begin{proof}
	1. \ Let $\left\{\pathp_i\right\}_{i \in I}$ be any set of paths in $\trt$ and let $\patha := \bigcap_{i \in I} \pathp_i \neq \emptyset$.  Pick any path $\pathp_0 \in \left\{\pathp_i\right\}_{i \in I}$.  That $\patha$ is a stem follows from the fact that $\patha \subseteq \pathp_0$ along with the fact that each path $\pathp_i$ is downward-closed.  To see that $\patha$ is a trunk it suffices to show, for all nodes $\ndu$ and $\ndv$ with $\ndu \in \patha$ and $\ndv \in \bridge{\ndu}$, that $\ndv \in \patha$.  Let $\ndu$ and $\ndv$ be any nodes for which $\ndu \in \patha$ and $\ndv \in \bridge{\ndu}$.  For each $i \in I$, since $\ndu \in \pathp_i$ then $\pathp_i \cap \bridge{\ndu} \neq \emptyset$ and so, since $\bridge{\ndu}$ is a bridge, $\pathp_i \cap \bridge{\ndu} = \bridge{\ndu}$, hence $v \in \pathp_i$.  It follows that $\ndv \in \patha$, as required.
	
	2. \ Let $\pathq$ be a path, $\left\{\pathp_i\right\}_{i \in I}$ be a set of paths in $\trt$, and $\pathb := \pathq \setminus \bigcup_{i \in I} \pathp_i \neq \emptyset$.  Let $\ndu \in \pathb$ and let $\ndv \in \pathq$ with $\ndu < \ndv$.  Suppose that $\ndv \not\in \pathb$.  Then $\ndv \in \pathp_j$ for some $j \in I$.  Then $\ndu \in \pathp_j$, as well, so that $\ndu \not\in \pathb$, a contradiction.  It follows that $\pathb$ is upward-closed in $\pathq$ hence $\pathb$ is a branch.
	
	To see that $\pathb$ is a limb it again suffices to show, for all nodes $\ndu$ and $\ndv$ with $\ndu \in \pathb$ and $\ndv \in \bridge{\ndu}$, that $\ndv \in \pathb$. So, let $\ndu$ and $\ndv$ be any nodes for which $\ndu \in \pathb$ and $\ndv \in \bridge{\ndu}$, but suppose, for a contradiction, that $\ndv \not\in \pathb$.  From $\ndv \in \bridge{\ndu}$,  
	it follows that $\ndv \in \pathp_j$ for some $j \in I$, hence $\bridge{\ndu} = \bridge{\ndv} \subseteq \pathp_j$. Then $\ndu \in \pathp_j$, so that $\ndu \not\in \pathb$, a contradiction.
	
	3. \ Let $\patha$ be a trunk, $\pathq$ a path that contains $\patha$, and $\pathb := \pathq \setminus \patha$.  We consider three cases (Cases 1 and 2 may overlap).
	
	Case 1: $\patha$ does not contain a greatest maximal bridge.  Let $\linacl$ be the set of maximal bridges in $\patha$.  Observe that for each $\segmj \in \linacl$ there exists a path $\pathp_{\segmj}$ and a node $\ndu_{\segmj} \in \patha$ such that $\segmj \subseteq \pathp_{\segmj}$ while $\ndu_{\segmj} \not\in \pathp_{\segmj}$.  Then $\pathb = \pathq \setminus \bigcup_{\segmj \in \linacl} \pathp_{\segmj}$.
	
	Case 2: $\pathb$ does not contain a least maximal bridge.  Let $\linbcl$ be the set of maximal bridges in $\pathb$.  Observe that for each $\segmj \in \linbcl$ there exists a path $\pathp_{\segmj}$ such that $\pathp_{\segmj} \cap \segmj = \emptyset$ while $\pathp_{\segmj} \cap \pathb \neq \emptyset$.  Then $\patha = \bigcap_{\segm \in \linbcl} \pathp_{\segm}$.
	
	Case 3: $\patha$ contains a greatest maximal bridge $\segm$ and $\pathb$ contains a least maximal bridge $\segmj$.  Since $\segm$ and $\segmj$ are distinct maximal bridges, there must exist a path $\pathp$ that contains $\segm$ but such that $\pathp \cap \segmj = \emptyset$.  Then $\pathp \cap \pathq = \patha$, as required.
\end{proof}

Trunks that are of the form in Part 1 of Theorem \ref{Thm:TrunksForm} will be called \defstyle{type I trunks}, and trunks that are of the form in Part 2 of the theorem will be called \defstyle{type II trunks}.  A limb will be called a \defstyle{type I limb} or \defstyle{type II limb} according to whether the trunk that complements it is a type I or type II trunk.  The set of paths $\left\{\pathp_i\right\}_{i \in I}$ from Theorem \ref{Thm:TrunksForm} will be said to \defstyle{generate} the corresponding trunk or limb as well as its complementary limb/trunk.  
A trunk $\patha$ will be called \defstyle{finitely generated} when it can be written in the form $\patha = \pathp \cap \pathq$ for some paths $\pathp$ and $\pathq$, and a limb $\pathb$ will be called \defstyle{finitely generated} when it can be written in the form $\pathb = \pathq \setminus \pathp$ for some paths $\pathq$ and $\pathp$.  If $\patha$ is a trunk that is contained in a path $\pathp$ then $\patha$ is finitely generated if and only if the limb $\pathp \setminus \patha$ is finitely generated.  Finitely generated trunks and limbs are both type I and type II but where the set of paths $\{\pathp_i\}_{i \in I}$ is finite.

\begin{lemma} \label{Thm:Trunk}
	Let $\patha$ be a finitely generated trunk in a tree $\trt$ and let $\pathp$ be a path such that $\patha \subseteq \pathp$.  Then there exists a path $\pathq$ such that $\patha = \pathp \cap \pathq$.
\end{lemma}

\begin{proof}
	Since $\patha$ is a finitely generated trunk, there exist paths $\pathb$ and $\pathc$ such that $\patha = \pathb \cap \pathc$.
		
	If $\pathp \cap \pathb = \patha$ then we are done, so consider the case where $\pathp \cap \pathb \supsetneq \patha$.  Let $\ndu \in (\pathp \cap \pathb) \setminus \patha$.  Suppose, for a contradiction, that $\pathp \cap \pathc \supsetneq \patha$.  Let $\ndv \in (\pathp \cap \pathc) \setminus \patha$.  Then $\min\left\{\ndu,\ndv\right\} \in (\pathb \cap \pathc) \setminus \patha$, a contradiction.  Hence $\pathp \cap \pathc = \patha$, as required.
\end{proof}

The following result characterises the notions of completeness \cmp{1} to \cmp{6}
 in terms of properties involving stems, branches, trunks, and limbs.

\begin{theorem} \label{Thm:CharacterisationOfCompleteness}
	Let $\trt$ be any tree. Then each of the following claims holds. 
	\begin{enumerate}
		\item
			$\trt$ is complete \cmp{1} if and only if every stem in $\trt$ has a supremum.
		\item
			$\trt$ is pathwise complete \cmp{2} if and only if every branch in $\trt$ has an infimum.
		\item
			$\trt$ is antichain complete \cmp{3} if and only if every type I trunk in $\trt$ has a greatest node.
		\item
			$\trt$ is branching complete \cmp{4} if and only if every finitely generated trunk in $\trt$ contains a greatest node.
		\item
			$\trt$ is weakly branching complete \cmp{5} if and only if every finitely generated limb in $\trt$ has an infimum.
		\item
			$\trt$ is branching complete \cmp{4} and weakly branching point complete \cmp{6} if and only if every trunk in $\trt$ has a supremum.
		\item
			$\trt$ is weakly branching complete \cmp{5} and weakly branching point complete \cmp{6} if and only if every limb in $\trt$ has an infimum.
	\end{enumerate}
\end{theorem}

\begin{proof}
	Claim 1: ($\Rightarrow$) \ Given a stem $\patha$ in $\trt$, the infimum of $T^{\geqslant \patha}$ can easily be verified to be the supremum of $\patha$.
	
	($\Leftarrow$) \ Given $X \subseteq T$ that is non-empty and bounded below, $\pathb := T^{\leqslant X}$ will be a stem in $\trt$.  Its supremum $\ndu$ can readily be verified to be the infimum of $X$.
	
	Claim 2: straightforward. 
	
	Claim 3: ($\Rightarrow$) \ Let $\trt$ be antichain complete and let $\patha := \bigcap_{i \in I} \pathp_i$, where $\{\pathp_i\}_{i \in I}$ is a set of paths, be a type I trunk.  Let $X$ be a maximal antichain in $\left(\bigcup_{i \in I} \pathp_i\right) \setminus \patha$.  Then $X$ is bounded below (by $\patha$) hence $X$ has an infimum $\ndu$.  Then $\ndu$ is the greatest node of $\patha$.
	
	($\Leftarrow$) \ Let $X$ be a non-empty antichain that is bounded below.  For each $\ndt \in X$, let $\pathp_{\ndt}$ be a path that passes through $\ndt$.  Then $\patha := \bigcap_{\ndt \in X} \pathp_{\ndt}$ is a type I trunk, hence $\patha$ has a greatest node $\ndu$.  This node $\ndu$ will be the required infimum of $X$.

	Claim 4: straightforward. 	
	
	Claim 5: straightforward (using Lemma \ref{Thm:Trunk}). 
	
	Claim 6: ($\Rightarrow$) \ Assume that $\trt$ satisfies \cmp{4} and \cmp{6}, and therefore also \cmp{7}.  Let $\patha$ be a trunk in $\trt$.  First, consider the case where $\patha$ is a type I trunk, say $\patha := \bigcap_{i \in I} \pathp_i$ for some set of paths $\{\pathp_i\}_{i \in I}$.  Fix $j \in I$ and for each $i \in I \setminus \{j\}$, let $\ndt_i$ be the greatest node in $\pathp_j \cap \pathp_i$, which exists by \cmp{4}.  By \cmp{7} the set $\{\ndt_i\}_{i \in I \setminus \{j\}}$ has an infimum $\ndu$, which is readily seen to be the supremum of $\patha$.  In the case where,  instead, $\patha$ is a type II trunk, it can similarly be shown that the supremum of the set of branching points that are formed on $\patha$ by the paths that generate $\patha$ is also a supremum of $\patha$ itself.
	
	($\Leftarrow$) \ Suppose that every trunk in $\trt$ has a supremum.  To see that \cmp{4} holds, let $\ndu$ and $\ndv$ be incomparable nodes and let $\patha$ and $\pathb$ be paths that contain $\ndu$ and $\ndv$.  The supremum of the trunk $\patha \cap \pathb$ will also be the infimum of $\{\ndu,\ndv\}$.  To see that \cmp{7}, and therefore also \cmp{6}, holds, let $X$ be a non-empty chain of branching points in $\trt$ and let $\pathp$ be a path for which $X \subseteq \pathp$.  For each $\ndt \in X$, let $\pathp_{\ndt}$ be a path for which $\ndt$ is the supremum of the trunk $\pathp \cap \pathp_{\ndt}$.  If $X$ is bounded below, then the supremum of the type I trunk $\bigcap_{\ndt \in X} \pathp_{\ndt}$ will be the infimum of $X$, while if $X$ is bounded above, then the supremum of the type II trunk $\pathp \setminus \bigcup_{\ndt \in X} \pathp_{\ndt}$ will be the supremum of $X$.
	
	Claim 7: The proof is similar to that of Claim 6.
\end{proof}

\begin{proposition} \label{Thm:PathwiseDedekindComplete}
	A tree $\trt$ is pathwise complete \cmp{2} if and only if $\trt$ is weakly branching complete \cmp{5} and weakly branching point complete \cmp{6} and each of the maximal bridges in $\trt$ is complete.
\end{proposition}

\begin{proof}
	($\Rightarrow$) \ Immediate.
	
	($\Leftarrow$) \ Let $\pathp$ be any path in $\trt$ and let $X \subseteq \pathp$ be non-empty and bounded below.
	
	Case 1: there exists a limb $\patha$ in $\trt$ such that $X$ is co-initial in $\patha$.  Then the infimum of $\patha$, which exists by Part 7 of Theorem \ref{Thm:CharacterisationOfCompleteness}, is also the infimum of $X$.
	
	Case 2: there exists a least maximal bridge $\pathb$ in $\trt$ such that $\pathb \cap X \neq \emptyset$, and a node $\ndb \in \pathb$ such that $b < x$ for each $x \in X$.  Since $\pathb$ is complete, the set $\pathb \cap X$ has an infimum which is also the infimum of $X$.
\end{proof}

\begin{corollary} \label{Thm:CompletenessCondensedTrees}
	Let $\trt$ be a condensed tree.
	\begin{enumerate}
		\item
			$\trt$ is pathwise complete \cmp{2} if and only if it is weakly branching complete \cmp{5} and weakly branching point complete \cmp{6}.
		\item
			$\trt$ is complete \cmp{1} if and only if it is branching complete \cmp{4} and weakly branching point complete \cmp{6}.
	\end{enumerate}
\end{corollary}

\begin{proof}
	1. This follows from Proposition \ref{Thm:PathwiseDedekindComplete} and the fact that every maximal bridge in a condensed tree is a singleton and hence complete.
	
	2. This follows from Theorem \ref{Thm:Dependences}, Theorem \ref{Thm:DedekindCompleteness} and Part 1 above.
\end{proof}

In what follows we need the notion of $<$-connected components in a forest.

\begin{definition}\label{def:schmerl-comp}
A \defstyle{$<$-connected component}\footnote{This definition comes from \cite{Schmerl}.}
(briefly, \defstyle{$<$-component}) of the forest $\forestf = \left(F;<\right)$  is a non-empty 
subset $\schm$ of $F$ such that: 

(1) if $\ndt \in \schm$, $\ndt' \leqslant \ndt$ and $\ndt' \leqslant \ndu$, then $\ndu \in \schm$; 

(2) $\schm$ is minimal (by inclusion) for the property (1). 
\end{definition}
It can be verified that the $<$-components of a forest coincide with the maximal subtrees of that forest.  The following characterisation will be needed when constructing completions of trees in Section \ref{sec:tree-completions}. 

\begin{proposition} \label{Thm:CompletenessComponents}
	Let $\trt$ be any tree.
	\begin{enumerate}
		\item
			$\trt$ is pathwise complete \cmp{2} if and only if for every stem $\pathx$ in $\trt$, each $<$-component of $\trt^{\geqslant \pathx}$ has a root. 
		\item
			$\trt$ is weakly branching complete \cmp{5} if and only if for every finitely generated trunk $\pathx$ in $\trt$, each $<$-component of $\trt^{\geqslant \pathx}$ has a root.
	\end{enumerate}
\end{proposition}

\begin{proof}
	1. ($\Rightarrow$) \ Let $\pathx$ be a stem in $\trt$.  If $\pathx$ has a greatest node $\ndu$ then $\ndu$ is the root of $\trt^{\geqslant \pathx}$, so let us consider the case where $\pathx$ has no greatest node.  Let $\tree'$ be any $<$-component of $\trt^{\geqslant \pathx}$ and let $\patha$ be a path in $\tree'$.  Then $\patha$ is a branch in $\trt$ for which $T^{< \patha} = \pathx$.  By Part 2 of Theorem~\ref{Thm:CharacterisationOfCompleteness}, $\patha$ has an infimum $\ndv$ in $\trt$.  Since $\pathx$ has no greatest node then $\ndv$ must be the least element of $\patha$.  It follows that $\ndv$ is the root of $\tree'$.
	
	($\Leftarrow$) \ By Part 2 of Theorem \ref{Thm:CharacterisationOfCompleteness}, it suffices to show that each branch in $\trt$ has an infimum.  Let $\patha$ be such a branch.  If $\patha$ has a least node then this least node will be its infimum, so let us consider the case where $\patha$ has no least node.  Let $\pathb := T^{< \patha}$ and let $\tree'$ be that $<$-component in $\trt^{\geqslant \pathb}$ that contains $\patha$.  By assumption, $\tree'$ has a root $\ndu$.  Since $\patha$ has no least node, $\ndu$ must be the greatest node of $\pathb$, from which it follows that $\ndu$ is the infimum of $\patha$.
	
	2. ($\Rightarrow$) \ By considering a finitely generated trunk $\pathx$, it can be shown, in a similar manner as in the forward direction of Part 1 of this proof, but using Part 5 rather than Part 2 of Theorem \ref{Thm:CharacterisationOfCompleteness}, that each $<$-component of $\trt^{\geqslant \pathx}$ has a root.
	
	($\Leftarrow$) \ This direction can be proved similarly to the backward direction of Part 1 of this proof, but again using Part 5, rather than Part 2, of Theorem \ref{Thm:CharacterisationOfCompleteness}, and considering a finitely generated limb $\patha$ rather than a branch, from which $\pathb := T^{<\patha}$ will be a finitely generated trunk, rather than a stem.
\end{proof}

\section{Completions of trees}
\label{sec:tree-completions}

We now turn to the problem of constructing completions of trees.  Intuitively, this will involve adding the necessary ``missing'' nodes to a tree so that it becomes complete,  in one sense or another.  The completions that will be considered here, are Dedekind completions, antichain completions, branching completions, pathwise Dedekind completions, weakly branching completions, and $\alpha$-fillings.  Theorem \ref{Thm:CharacterisationOfCompleteness} and Proposition \ref{Thm:CompletenessComponents} will be key in constructing these completions in intuitive terms.

\subsection{Dedekind completions}

We start by outlining a construction for obtaining the Dedekind completion of a tree that is given in \cite{Warren} (applied there to the wider class of so called `cycle-free partial orders').
 Given a tree $\trt$, a nonempty subset $X$ of $T$ is called a \defstyle{Dedekind ideal} when $X$ is bounded above and $L(U(X)) = X$; $X$ is called a \defstyle{principal ideal} when it has the form $X = T^{\leqslant y}$ for some $y \in T$.  Let $\mathcal{I}\left(\trt\right)$ denote the set that consists of all Dedekind ideals of $\trt$ and let us define $\trt^{\mathcal{I}} := \left(\mathcal{I}\left(\trt\right);\subseteq\right)$.

\begin{fact}\cite[Lemma 2.2.5]{Warren}
	For any tree $\trt$, the structure $\trt^{\mathcal{I}}$ is a complete tree.\footnote{In \cite{Warren}, Dedekind completeness is defined as the property that each ideal is principal.  This amounts to the same as the definition of (Dedekind) completeness that is used in this paper.}
\end{fact}

The tree $\trt^{\mathcal{I}}$ can be related to $\trt$ as follows: define $f : \trt \rightarrow \trt^{\mathcal{I}}$ by $x \mapsto T^{\leqslant x}$, i.e.,~each node in $\trt$ is mapped to the principal ideal that it generates.

\begin{fact}\cite[Lemmas 2.2.3, 2.2.5, 2.2.7]{Warren} 
\label{Thm:DedekindCompletionMinimal}
The mapping 
	$f$ is an isomorphic embedding of $\trt$ into $\trt^{\mathcal{I}}$, and $\trt^{\mathcal{I}}$ is a substructure of every complete tree that extends $f\left[\trt\right]$.
\end{fact}

We note that the above construction is an adaptation of the usual \emph{Dede\-kind-MacNeille completion} (see e.g.,~\cite{MacNeille}, \cite{BurrisSankappanavar}) of a partially ordered set, but with the following differences: the Dedekind-MacNeille completion of a tree $\trt$ takes for its underlying set \textit{all} subsets of $T$ for which $L(U(X)) = X$, not only those sets $X$ that are non-empty and bounded above, and the resulting structure is (up to isomorphism) the smallest complete \textit{lattice} that contains $\trt$.

We now describe how to construct the Dedekind
completion of a tree in more intuitive terms.\footnote{A similar construction is given in  \cite{Droste85}, \S5.}   Given a tree $\trt := \left(T;<\right)$, let $\mathcal{S}$ denote the set that consists of all stems in $\trt$ that do not have a supremum.  For each $\stem \in \mathcal{S}$, let $\ndt_{\stem}$ denote a new node that is not in $\trt$, and let $T^{DC} := T \cup \left\{\ndt_{\stem} : \stem \in \mathcal{S}\right\}$.  Let $<^{DC}$ be the relation
\begin{multline} \label{Eqn:CompletionOrdering}
	\left\{(x,y) \in T \times T: x < y\right\} \cup \left\{\left(\ndt_{\stemr},\ndt_{\stem}\right) : \stemr,\stem \in \mathcal{S} \ \text{with} \ \stemr \subsetneq \stem \right\} \cup \\
	\bigcup_{\stem \in \mathcal{S}}\left\{\left(\ndu,\ndt_{\stem}\right) : u \in \stem \right\} \cup \bigcup_{\stem \in \mathcal{S}}\left\{\left(\ndt_{\stem},\ndu\right) : \ndu \ \text{is any upper bound of} \ \stem \right\}
\end{multline}
on $T^{DC}$ and let us denote the structure $\left(T^{DC};<^{DC}\right)$ as $\trt^{DC}$.

\begin{proposition} \label{Thm:DedekindCompletion}
	The structure $\trt^{DC}$ is a complete tree that contains $\trt$.
\end{proposition}

\begin{proof}
	It is straightforward to check that $<^{DC}$ is irreflexive, transitive, downward-linear, and downward-connected, hence $\trt^{DC}$ is a tree.
	
	By Theorem \ref{Thm:CharacterisationOfCompleteness}, we have only to show that every stem in $\trt^{DC}$ has a supremum.  Observe firstly that every stem $\stem$ that is also a stem in $\trt$, has a supremum in $\trt^{DC}$, namely either its supremum $\sup(\stem)$ in $\trt$, or $\ndt_{\stem}$.
	
	Next, let $\stemr$ be any stem in $\trt^{DC}$ that is not a stem in $\trt$.  If $\stemr$ has a greatest element, there is nothing to prove.  Otherwise, consider the set $\stemr^- := \stemr \cap T$, which is a stem in $\trt$.  It suffices to show that $\stemr^-$ is cofinal in $\stemr$, that is, for every $\ndu \in \stemr$, there is a $\ndt \in \stemr^-$ such that $\ndu \leqslant^{DC} \ndt$; then $\ndt_{\stemr^-}$, the supremum of $\stemr^-$ in $\trt^{DC}$, will also be the supremum of $\stemr$ in $\trt^{DC}$.
	
	Indeed, suppose that $\ndu$ is an element of $\stemr$ such that $\ndt <^{DC} \ndu$ for all $\ndt \in \stemr^-$.  Then $\ndu \not\in \stemr^-$, hence $\ndu = \ndt_{\stem}$ for some stem $\stem$ in $\trt$.  Since $\stem \subseteq \stemr$ then $\stem \subseteq \stemr^-$, and it is easily verified that $\stemr^- \subseteq \stem$, hence $\stem = \stemr^-$, so that $\ndu = \ndt_{\stem} = \ndt_{\stemr^-}$.  Now consider any node $\ndu'$ in $\stemr$ for which $\ndu <^{DC} \ndu'$; such $\ndu'$ exists because $\stemr$ has no greatest element.  Then $\ndt <^{DC} \ndu'$ for all $\ndt \in \stemr^-$ and it again follows that $\ndu' = \ndt_{\stemr^-}$ hence $\ndu' = \ndu$, which contradicts the irreflexivity of $<^{DC}$.
\end{proof}

The tree $\trt^{DC}$ will be called the \defstyle{Dedekind completion} of $\trt$.  By Theorem \ref{Thm:CharacterisationOfCompleteness}, every complete tree must have a supremum for each of its stems.  It follows that $\trt^{DC}$ can be embedded in every complete tree that extends $\trt$.  Using Fact \ref{Thm:DedekindCompletionMinimal}, it therefore must be the case that $\trt^{DC} \cong \trt^{\mathcal{I}}$.

\subsection{Antichain completions and branching completions}

Given a tree $\trt := \left(T;<\right)$ with Dedekind completion $\trt^{DC}$, the antichain completion and branching completion\footnote{\cite{Barham} defines branching completions, there called \textit{ramification completions}, in this way.} of $\trt$ can be obtained respectively as
$$\bigcap \left\{ \trs : \trs \ \text{is an antichain complete tree such that} \ \trt \subseteq \trs \subseteq \trt^{DC} \right\}$$
and
$$\bigcap \left\{ \trs : \trs \ \text{is a branching complete tree such that} \ \trt \subseteq \trs \subseteq \trt^{DC} \right\}.$$
Theorem \ref{Thm:CharacterisationOfCompleteness} again suggests intuitively simple constructions for these completions, this time in terms of the trunks of $\trt$.

Let $\mathcal{T}_{\mathrm{I}}$ be the set of all type I trunks in $\trt$ that do not contain a greatest node, and let $\mathcal{T}_{\mathrm{fin}}$ be the set of all finitely generated trunks in $\trt$ that do not contain a greatest node.  For each trunk $\stem$ in $\trt$, let $t_{\stem}$ denote a new node that is not in $T$, and let us define
\begin{eqnarray*}
	T^{AC} & := & T \cup \left\{t_{\stem} : \stem \in \mathcal{T}_{\mathrm{I}} \right\},  \\ 
	T^{BC} & := & T \cup \left\{t_{\stem} : \stem \in \mathcal{T}_{\mathrm{fin}} \right\}.
\end{eqnarray*}
Now, we define the relations $<^{AC}$ and $<^{BC}$ as in (\ref{Eqn:CompletionOrdering}) but respectively with the sets $\mathcal{T}_{\mathrm{I}}$ and $\mathcal{T}_{\mathrm{fin}}$ replacing the set $\mathcal{S}$. Finally, let $\trt^{AC} := \left(T^{AC};<^{AC}\right)$ and $\trt^{BC} := \left(T^{BC};<^{BC}\right)$.

\begin{proposition} \label{Thm:AntichainBranchingCompletion}
	Let $\trt$ be any tree.
	\begin{enumerate}
		\item
			$\trt^{AC}$ is an antichain complete tree that contains $\trt$.
		\item
			$\trt^{BC}$ is a branching complete tree that contains $\trt$.
	\end{enumerate}
\end{proposition}

\begin{proof}
	1. To show that $\trt^{AC}$ is a tree, we can use an argument identical to the one  used in the proof of Proposition \ref{Thm:DedekindCompletion} to show that $\trt^{DC}$ was a tree.
		
	Let $\stemr$ be a type I trunk in $\trt^{AC}$, say $\stemr = \bigcap_{i \in I} \pathp_i$ for some set of paths $\left\{\pathp_i\right\}_{i \in I}$ in $\trt^{AC}$.  Let $\stemr^- := \stemr \cap T$ and $\pathp_i^- := \pathp_i \cap T$ for each $i \in I$.  Then for each $i$, $\pathp_i^-$ is a path in $\trt$, and $\stemr^- = \bigcap_{i \in I} \pathp_i^-$ is a type I trunk in $\trt$.
	
	First consider the case where $\stemr^-$ has a greatest node $\ndu$.  Then $\ndu$ is also the greatest node of $\stemr$.  To see this, suppose, for a contradiction, that there exists $\ndv \in \stemr$ for which $\ndu <^{AC} \ndv$.  Then $\ndv \not\in T$ (else $\ndv \in \stemr^-$, a contradiction) hence $\ndv = \ndt_{\stem}$ for some trunk $\stem$ in $\mathcal{T}_{\mathrm{I}}$.  Since also $\stem \subseteq \stemr$ (because $\ndv \in \stemr$) then $\stem \subseteq \stemr^-$, and since $\ndu <^{AC} \ndt_{\stem}$ then $\stemr^- \subseteq \stem$, hence $\stem = \stemr^-$.  This contradicts the fact that $\stem$ does not have a greatest node.
	
	Next, consider the case where $\stemr^-$ does not have a greatest node.  Then $T^{AC}$ contains the node $\ndt_{\stemr^-}$, and $\ndu := \ndt_{\stemr^-}$ will be the greatest node of $\stemr$.  To see this, first note, from the way that $<^{AC}$ is defined, that $u \in \pathp_i$ for each $i \in I$, hence $u \in \stemr$.  Next, suppose  again, for a contradiction, that there exists $\ndv \in \stemr$ for which $\ndu <^{AC} \ndv$.  It follows again that $\ndv \not\in T$, for if $\ndv \in T$ then, as above, $\ndv \in \stemr^-$ hence $\ndv <^{AC} \ndt_{\stemr^-}$, which contradicts the fact that $\ndu <^{AC} \ndv$.  Hence, $\ndv = \ndt_{\stem}$ for some trunk $\stem$ in $\mathcal{T}_{\mathrm{I}}$.  Again, it follows from $\ndv \in \stemr$ that $\stem \subseteq \stemr$, hence $\stem \subseteq \stemr^-$, while from $\ndt_{\stemr^-} = \ndu <^{AC} \ndv = \ndt_{\stem}$ can be concluded that $\stemr^- \subsetneq \stem$, a contradiction.
	
	The proof of Claim 2 is similar. 
\end{proof}

The trees $\trt^{AC}$ and $\trt^{BC}$ will be called, respectively,  the \defstyle{antichain completion} and \defstyle{branching completion} of $\trt$.  Using Theorem \ref{Thm:CharacterisationOfCompleteness}, it follows that $\trt^{AC}$ can be embedded in every antichain complete tree that extends $\trt$, and $\trt^{BC}$ can be embedded in every branching complete tree that extends $\trt$.

\subsection{Pathwise Dedekind completions and weakly branching completions}

The case of constructing a pathwise Dedekind completion of a tree poses a minor complication, in that the construction is not deterministic.  In general, there need not be a unique minimal pathwise complete tree that extends a given tree.  In terms of Theorem \ref{Thm:CharacterisationOfCompleteness}, to obtain a pathwise complete tree from a given tree it suffices to add either missing least nodes to the tree's branches, or missing greatest nodes to its stems, and there is some freedom in how one can go about doing this.

Instead, we will use Proposition \ref{Thm:CompletenessComponents} as the basis for our construction.  The reason for performing the construction using $<$-components,  rather than using branches as suggested by Theorem \ref{Thm:CharacterisationOfCompleteness}, is that, if one were to simply add nodes as infima to branches that are without infima, it could happen that multiple nodes get added to the same branch since different branches may have the same set of lower bounds.

The construction presented below will produce, from any tree $\trt$, a minimal pathwise complete tree $\trt^{PDC}$, without introducing any additional bran\-ch\-ing points. Thus, $\trt^{PDC}$ will be the most conservative pathwise Dedekind completion of $\trt$ in the sense that, amongst all pathwise Dedekind completions of $\trt$, $\trt^{PDC}$ will be the one that is furthest away from being branching complete.

Similar observations hold when producing weakly branching completions from a tree, where, again, a tree need not have a unique minimal weakly branching completion.  Here we will again base the construction on Proposition \ref{Thm:CompletenessComponents}, rather than on Theorem \ref{Thm:CharacterisationOfCompleteness}, to obtain from a tree $\trt$ a minimal weakly branching complete tree $\trt^{WBC}$ that does not contain any new branching points and which, amongst all weakly branching complete trees that extend $\trt$, will be the one that is furthest away from being branching complete.

\smallskip

Let $\trt := \left(T;<\right)$ be any tree.  Let $\mathcal{S}'$ be the set of all stems in $\trt$, and let
\begin{equation} \label{Eqn:ClassC}
	\mathcal{C}(\mathcal{S}') := \bigcup_{\stem \in \mathcal{S}'}\left\{ \schm : \schm \ \text{is a $<$-component without root in} \ \trt^{\geqslant \stem} \right\}.
\end{equation}
For each $\schm \in \mathcal{C}(\mathcal{S}')$, let $\ndt_{\schm}$ denote a new node that is not in $\trt$, and define $T^{PDC} := T \cup \left\{\ndt_{\schm} : \schm \in \mathcal{C}(\mathcal{S}')\right\}$.  Let $<^{PDC}$ be the following relation on $T^{PDC}$:  
\begin{multline}
	\left\{(x,y) \in T \times T: x < y\right\} \cup \left\{\left(\ndt_{\schm},\ndt_{\schmd}\right) : \schm, \schmd \in \mathcal{C}(\mathcal{S}') \ \text{with} \ \schmd \subsetneq \schm \right\} \cup \\
	\bigcup_{\schm \in \mathcal{C}(\mathcal{S}')}\left\{\left(\ndu,\ndt_{\schm}\right) : \ndu \ \text{is any lower bound of} \ \schm \right\} \cup \bigcup_{\schm \in \mathcal{C}(\mathcal{S}')}\left\{\left(\ndt_{\schm},\ndu\right) : \ndu \in \schm \right\}.  \label{Eqn:RelationPDC}
\end{multline}
Now, we define $\trt^{PDC} := \left(T^{PDC};<^{PDC}\right)$.

\medskip

The tree $\trt^{WBC}:= \left(T^{WBC};<^{WBC}\right)$ is defined similarly, as follows.  Let $\mathcal{T}'$ be the set of all finitely generated trunks in $\trt$, and define $\mathcal{C}(\mathcal{T}')$ as in (\ref{Eqn:ClassC}) but using the set $\mathcal{T}'$ in place of $\mathcal{S}'$.  For each $\schm \in \mathcal{C}(\mathcal{T}')$, let $\ndt_{\schm}$ again denote a new node that is not in $\trt$, and again define $T^{WBC} := T \cup \left\{\ndt_{\schm} : \schm \in \mathcal{C}(\mathcal{T}')\right\}$.  Finally, we define the relation $<^{WBC}$ as in (\ref{Eqn:RelationPDC}), but again using the set $\mathcal{T}'$ instead of $\mathcal{S}'$.

\begin{proposition} \label{Thm:PathwiseDedekindCompletion}
	Let $\trt$ be any tree.
	\begin{enumerate}
		\item
			$\trt^{PDC}$ is a pathwise complete tree that contains $\trt$.
		\item
			$\trt^{WBC}$ is a weakly branching complete tree that contains $\trt$.
	\end{enumerate}
\end{proposition}

\begin{proof}
	1. The verification that $\trt^{PDC}$ is a tree is straightforward.
	
	Let $\stemr$ be a stem in $\trt^{PDC}$ and let $\schm$ be a $<^{PDC}$-component in the forest $\left(\trt^{PDC}\right)^{\geqslant^{PDC}\,  \stemr}$.  By Proposition \ref{Thm:CompletenessComponents}, it suffices to show that $\schm$ has a root.  Let $\schm^- := \schm \cap T$.  Two cases will be considered.
	
	Case 1. Suppose that $\schm^-$ has a root $\ndu$.  Then $\ndu$ is also a root in $\schm$.  To see this, it suffices to show that $\ndu <^{PDC} \ndt_{\schmd}$ for each node in $\schm$ of the form $\ndt_{\schmd}$, since it is immediate that $\ndu <^{PDC} \ndv$ for each $\ndv \in \schm^-$.  Suppose, for a contradiction, that $\ndu \not\smile^{PDC} \ndt_{\schmd}$ or $\ndt_{\schmd} <^{PDC} \ndu$.
	If $\ndu \not\smile^{PDC} \ndt_{\schmd}$ then $\ndu \not<^{PDC} \schmd$ and $\ndu \not\in \schmd$.  It follows that $\ndu \not\smile^{PDC} \ndv$ for each $\ndv \in \schmd$, which contradicts the fact that $\schmd \subseteq \schm^-$ while $\ndu$ is the root of $\schm^-$.
	On the other hand, if $\ndt_{\schmd} <^{PDC} \ndu$ then $u \in \schmd$.  Since $\ndu$ is the root of $\schm^-$ then $\schm^- \subseteq \schmd$, and since $\ndt_{\schmd} \in \schm$ then $\schmd \subseteq \schm$ hence $\schmd \subseteq \schm^-$, so that $\schmd = \schm^-$.  Hence $\ndt_{\schmd} = \ndt_{\schm^-}$, but this contradicts the assumption that $\schm^-$ has a root.
	
	Case 2. Suppose that $\schm^-$ does not have a root.  Then the node $\ndt_{\schm^-}$ is the root of $\schm$.  That $\ndt_{\schm^-} \in \schm$ follows from $\stemr \leqslant^{PDC} \ \ndt_{\schm^-} <^{PDC} \schm^-$, along with the fact that $\schm \supseteq \schm^-$ and $\schm$ is a $<^{PDC}$-component in $\left(\trt^{PDC}\right)^{\geqslant^{PDC} \, \stemr}$.  It now suffices to show that $\ndt_{\schm^-} <^{PDC} \ndt_{\schmd}$ for each $\ndt_{\schmd} \in \schm$ with $\ndt_{\schmd} \neq \ndt_{\schm^-}$.  Indeed, if $\ndt_{\schmd} \in \schm$ then $\schmd \subseteq \schm$ hence $\schmd \subseteq \schm^-$, and since, in addition, $\ndt_{\schmd} \neq \ndt_{\schm^-}$, then $\schmd \subsetneq \schm^-$ from which $\ndt_{\schm^-} <^{PDC} \ndt_{\schmd}$.
	
	2. The proof is identical to that of Claim 1, except for using finitely generated trunks, rather than stems.
\end{proof}

Using Proposition \ref{Thm:CompletenessComponents}, it follows that $\trt^{PDC}$ is minimal in the class of pathwise complete trees that extend $\trt$, and $\trt^{PDC}$ can be embedded in every pathwise complete tree that extends $\trt$ and has the same branching points as $\trt$.  Similarly, $\trt^{WBC}$ is minimal in the class of weakly branching complete trees that extend $\trt$, and $\trt^{WBC}$ can be embedded in every weakly branching complete tree that extends $\trt$ and has the same branching points as $\trt$.

\smallskip

Lastly, we raise the following question, which we leave open: for which families $\mathcal{F}$ of sets of nodes in a tree, does the generic notion of $\mathcal{F}$-completeness give rise to a generic construction of an $\mathcal{F}$-completion of the tree that can be obtained in a way similar to the constructions presented here?

\subsection{$\alpha$-fillings}

The Dedekind completion that was proposed above need not preserve some natural properties, such as denseness.  For example, consider a tree $\trt$ that consists of a copy of the rationals $\eta$, on top of which two incomparable copies of $1 + \eta$ are appended.  $\trt$ has two paths, each of order type $\eta + 1 + \eta \cong \eta$.  Both of these paths become copies of the linear order $\lambda + 2 + \lambda$ in the Dedekind completion of $\trt$.  The paths in $\trt$ are thus dense linear orders whereas the paths in the Dedekind completion of $\trt$ are not dense.

The question therefore arises how to modify the construction to produce a dense Dedekind completion?  Here we propose an alternative general construction that, in particular, does that.  The construction is non-deterministic, and is defined as follows.

Given a condensed tree $\trs := \left(S;<_{\trs}\right)$, a set of linear orders $\linacl := \left\{\lina_i\right\}_{i \in I}$, and a function $f : S \rightarrow I$, the \defstyle{$f$-product} $\trs \times_f \linacl := \left( \lvert \trs \times_f \linacl \rvert; < \right)$ is the structure that is defined as follows:
\begin{itemize}
	\item
		$\displaystyle \lvert \trs \times_f \linacl \rvert := \bigcup_{\ndt \in S} \left(\left\{\ndt\right\} \times \lvert\lina_{f(\ndt)}\rvert\right)$;
	\item
		for $\left(\ndt_1,x_1\right),\left(\ndt_2,x_2\right) \in \lvert \trs \times_f \linacl \rvert$,
		$$\left(\ndt_1,x_1\right) < \left(\ndt_2,x_2\right) \ \Longleftrightarrow \ \ndt_1 <_{\trs} \ndt_2 \ \text{or} \ \left(\ndt_1 = \ndt_2 \ \text{and} \ x_1 <_{\lina_{f\left(t_1\right)}} x_2\right).$$
\end{itemize}
Informally, the structure $\trs \times_f \linacl$ is obtained from $\trs$ by replacing each node $\ndt$ in $\trs$ with the linear order $\lina_{f\left(\ndt\right)}$.  It is easily seen that $\trs \times_f \linacl$ is a tree, that the condensation (cf. \cite{GorankoKellermanZanardo2021a})
of $\trs \times_f \linacl$ is isomorphic to $\trs$, and that the distinct maximal bridges in $\trs \times_f \linacl$ are, up to isomorphism, precisely the elements of $f\left(S\right) \subseteq \mathcal{A}$.  The tree $\trs \times_f \linacl$ is, strictly speaking, not an extension of $\trs$ itself, because the elements of $\trs \times_f \linacl$ are ordered pairs while the elements of $\trs$ are not.  However, $\trs \times_f \linacl$ \textit{is} isomorphic to a tree $\trs'$ that extends $\trs$, and $\trs \times_f \linacl$ will, for the sake of simplicity, be identified with this tree $\trs'$.

Now, let $\alpha$ be any dense linear order without endpoints and let $\trt = (T;<)$ be a condensed tree.  Take $\lina_0 := \alpha$ and $\lina_1 := \alpha+1$.  Define an \defstyle{$\alpha$-filling} of $\tree$ to be any tree of the form
\begin{equation} \label{Eqn:AlphaFilling}
	\trt \times_f \{\lina_0,\lina_1\}
\end{equation}
where $f : T \rightarrow \{0,1\}$ is any function.  Intuitively, an $\alpha$-filling of $\trt$ is obtained by replacing each node in $\trt$ with either $\alpha$ or $\alpha+1$.  Thus, the tree $\trt$ provides the branching structure, and $\alpha$ provides the filling pattern inside maximal bridges.  A tree that is of the form (\ref{Eqn:AlphaFilling}) will also be said to be a \defstyle{locally $\alpha$-tree}.  If $f(x) = \lina_0$ for each leaf in $\trt$, and $f(x) = \lina_1$ otherwise, then $\trt \times_f \{\lina_0,\lina_1\}$ will be called the \defstyle{full $\alpha$-filling} of $\trt$.

In the example of the tree $\trt$ described above, that consists of a copy of $\eta$ with two disjoint copies of $1+\eta$ appended on top of it, producing the full $\eta$-filling from the condensation of $\trt$, will have the effect of collapsing the two least elements of each of the copies of $1+\eta$ in $\trt$, into a single branching node.  The resulting $\eta$-filling will be both dense and branching complete.

A linear order is called \defstyle{continuous} when it is dense, without endpoints, and complete.  A linear order $\lina := \left(A;<\right)$ is called \defstyle{separable} if each subset $X$ of $A$ that is dense in $\lina$ (i.e.,~whenever $a,b \in A$ with $a < b$, there exists $x \in X$ such that $a < x < b$) is at most countable.  Recall the following two facts:
\begin{itemize}
	\item
		Cantor's Theorem (see e.g.~\cite[Theorem 2.8]{Rosenstein}): every countable dense linear order without endpoints is isomorphic to $\eta$; and
	\item
		a linear order is isomorphic to $\lambda$ if and only if it is continuous and separable (see e.g.~\cite[Theorem 2.30]{Rosenstein}).
\end{itemize}
The following two special cases of $\alpha$-fillings are now singled out:
\begin{itemize}
	\item
		When $\alpha = \eta$. Every full $\eta$-filling of a condensed tree $\trt$ has the property that each of its paths is a dense linear order without endpoints.  Moreover, if each path in $\trt$ is at most countable then, by Cantor's Theorem, each path in the full $\eta$-filling of $\trt$ will be isomorphic to $\eta$.
	\item
		When $\alpha = \lambda$. Given a condensed pathwise complete tree $\trt$, each path in the full $\lambda$-filling of $\trt$ will be a continuous linear order.  Moreover, if each path in $\trt$ is at most countable then, by the above characterisation of $\lambda$, it follows that each path in the full $\lambda$-filling of $\trt$ will be isomorphic to $\lambda$.  Finally, if $\trt$ is branching complete then, by Theorem \ref{Thm:DedekindCompleteness}, its full $\lambda$-filling will be complete.
\end{itemize}

\section{Concluding remarks}
\label{sec:concluding} 

This work continues the study of the general theory of trees initiated in 
 \cite{GorankoKellermanZanardo2021a}. 
Further planned work on this project will include: 

\begin{itemize}

\item the study of general operations on trees, such as sums and products, thus extending classical studies of ordinal arithmetic (Cantor, Sierpinski, and others) and, more generally, operations on linear orderings 
(cf. \cite{Rosenstein});  

\item the study of classes of trees generated by applying such operations, and their structural and logical theories.
\end{itemize}

The ultimate goal of this project is a systematic development of a structural theory of trees. One intended application of this study is to obtain structural characterisations of elementary equivalence and other logical equivalences of trees and to use them to obtain new axiomatisations and decidability or undecidability results for logical theories of important classes of trees.

\section*{Acknowledgements}

We thank the referee for the careful reading and helpful comments and suggestions on the paper.

\providecommand{\bysame}{\leavevmode\hbox to3em{\hrulefill}\thinspace}
\providecommand{\MR}{\relax\ifhmode\unskip\space\fi MR }
\providecommand{\MRhref}[2]{%
  \href{http://www.ams.org/mathscinet-getitem?mr=#1}{#2}
}
\providecommand{\href}[2]{#2}

\end{document}